\documentclass[uplatex,12pt]{article}% upLaTeX 
\usepackage{amsmath}
\usepackage{amssymb}
\usepackage{latexsym}
\usepackage[dvips]{color}
\usepackage{fancybox}
\usepackage{graphics}
\usepackage{amsthm}
\usepackage{tabularx}
\usepackage{amsthm}
\usepackage{subcaption}
\usepackage{mathrsfs} %花文字
\usepackage{colortbl}
\usepackage{tikz-cd}
\usepackage{bm}
\usepackage{hyperref}

\usepackage[a4paper]{geometry}
\usepackage{mathrsfs}
\usepackage{comment}
\usepackage{accsupp}

\def\diam{\mathop{\mathrm{diam}}}
\def\diag{\mathop{\mathrm{diag}}}

\def\div{\mathop{\mathrm{div}}}

\def\<{\mathop{\textless}}
\def\>{\mathop{\textgreater}}

\def\card{\mathop{\rm{card}}}

\theoremstyle{definition}
\newtheorem{thr}{Theorem}[section]

\newtheorem*{ex*}{Example}
\newtheorem{rem}[thr]{Remark}
\newtheorem{defi}[thr]{Definition}
\newtheorem{lem}[thr]{Lemma}

\newtheorem*{pf*}{Proof}

\newtheorem{assume}{Assumption}{\bf}{\it}
\newtheorem*{ass*}{Assumption}
\newtheorem{Cond}[thr]{Condition}

\newtheorem*{exer*}{Exercise}

\newtheorem*{memo*}{Memo}

\newtheorem*{pfth*}{Proof of Theorem 2.6}
\newtheorem*{pfth1*}{Proof of Theorem 2.7}

\newtheorem*{lem21*}{Lemma 2.1}
\newtheorem*{th15*}{Theorem 1.5}
\newtheorem*{th21*}{Theorem 2.1}
\newtheorem*{th25*}{Theorem 2.5}
\newtheorem*{th26*}{Theorem 2.6}

\usepackage{makeidx}
\makeindex

\allowdisplaybreaks[1]

\newcounter{sone}
\newcounter{stwo}
\newcounter{sthree}
\newcounter{sfour}
\newcounter{sfive}
\newcounter{ssix}
\newcounter{lone}
\newcounter{ltwo}
\newcounter{lthree}
\newcounter{lfour}
\newcounter{lfive}
\newcounter{lsix}
\newcounter{lseven}
\newcounter{leight}
\makeatletter
    
    \@addtoreset{equation}{section}
  \makeatother

\begin{document}

%%% 節番号 %%%
\setcounter{section}{0}
%%% titlefoot.sty 形式 %%%
\title{Special topics in FEMs \\ [1ex] \large Note on a weakly over-penalised symmetric interior penalty method on anisotropic meshes for the Poisson equation, ver. 1}
%\subtitle{Note on a weakly over-penalised symmetric interior penalty method on anisotropic meshes for the Poisson equation, ver. 1}
\author{Hiroki Ishizaka \thanks{Team FEM, Matsuyama, Japan,  h.ishizaka005@gmail.com}}
%\institute{Hiroki Ishizaka \at
    %          Team FEM, Matsuyama, Japan \\
      %        \href{h.ishizaka005@gmail.com} 
%}

\date{\today}

\maketitle
\pagestyle{plain}

\begin{abstract}
The purpose is to make an easy-to-understand note of "Special Topics in Finite Element Methods." There might be typos and mistakes. Therefore, I do not take any responsibility for unauthorised use.


\end{abstract}

\section{Introduction}
In this note, we investigate a weakly over-penalised symmetric interior penalty (WOPSIP) method for the Poisson equations on anisotropic meshes. Brenner et al. first proposed a WOPSIP method \cite{BreOweSun08}, and several further works have considered similar techniques \cite{BarBre14,Bre15,BreOweSun12}. Analysis of the WOPSIP method on anisotromeshes is studied in \cite{Ish24} for the Stokes equations.

\section{Preliminaries}
Throughout, we denote by $c$ a constant independent of $h$ (defined later) and of the angles and aspect ratios of simplices unless specified otherwise, {and all constants $c$ are bounded if the maximum angle is bounded}. These values may change in each context. The notation $\mathbb{R}_+$ denotes the set of positive real numbers.

\subsection{Continuous problems}
Let $\Omega \subset \mathbb{R}^d$, $d \in \{ 2 , 3 \}$ be a bounded polyhedral domain. Furthermore, we assume that $\Omega$ is convex if necessary. The Poisson problem is to find $u: \Omega \to \mathbb{R}$ such that
\begin{align}
\displaystyle
- \varDelta u  = f \quad \text{in $\Omega$}, \quad u = 0 \quad \text{on $\partial \Omega$}, \label{poisson_eq}
\end{align}
where $f \in L^2(\Omega)$ is a given function. The variational formulation for the Poisson equations \eqref{poisson_eq} is as follows. Find $u \in H_0^1(\Omega)$ such that
\begin{align}
\displaystyle
a(u,\varphi) := \int_{\Omega} \nabla u \cdot \nabla \varphi dx &= \int_{\Omega} f  \varphi dx \quad \forall \varphi \in H_0^1(\Omega). \label{poisson_weak}
\end{align}
By the Lax--Milgram lemma, there exists a unique solution $u \in H_0^1(\Omega)$ for any $f \in L^2(\Omega)$ and it holds that
\begin{align*}
\displaystyle
| u |_{H^1(\Omega)} \leq C_P(\Omega) \| f \|,
\end{align*}
where $C_P(\Omega)$ is the Poincar$\rm{\acute{e}}$ constant depending on $\Omega$. 
%\begin{comment}
Furthermore, if $\Omega$ is convex, then $u \in H^2(\Omega)$ and 
\begin{align}
\displaystyle
| u |_{H^2(\Omega)} \leq \| \varDelta u \|. \label{poisson=reg}
\end{align}
The proof can be found in, for example, \cite[Theorem 3.1.1.2, Theorem 3.2.1.2]{Gri11}. 

\subsection{Meshes, mesh faces, averages and jumps} \label{regularmesh}
 Let $\mathbb{T}_h = \{ T \}$ be a simplicial mesh of $\overline{\Omega}$ made up of closed $d$-simplices such as $\overline{\Omega} = \bigcup_{T \in \mathbb{T}_h} T$ with $h := \max_{T \in \mathbb{T}_h} h_{T}$, where $ h_{T} := \diam(T)$. For simplicity, we assume that $\mathbb{T}_h$ is conformal: that is, $\mathbb{T}_h$ is a simplicial mesh of $\overline{\Omega}$ without hanging nodes.

Let $\mathcal{F}_h^i$ be the set of interior faces and $\mathcal{F}_h^{\partial}$ the set of the faces on the boundary $\partial \Omega$. We set $\mathcal{F}_h := \mathcal{F}_h^i \cup \mathcal{F}_h^{\partial}$. For any $F \in \mathcal{F}_h$, we define the unit normal $n_F$ to $F$ as follows. (\roman{sone}) If  $F \in \mathcal{F}_h^i$ with $F = T_1 \cap T_2$, $T_1,T_2 \in \mathbb{T}_h$, let $n_1$ and $n_2$ be the outward unit normals of $T_1$ and $T_2$, respectively. Then, $n_F$ is either of $\{ n_1 , n_2\}$; (\roman{stwo}) If $F \in \mathcal{F}_h^{\partial}$, $n_F$ is the unit outward normal $n$ to $\partial \Omega$. For a simplex $T \subset \mathbb{R}^d$, let $\mathcal{F}_{T}$ be the collection of the faces of $T$.

Here, we consider $\mathbb{R}^q$-valued functions for some $q \in \mathbb{N}$. We define a broken (piecewise) Hilbert space as
\begin{align*}
\displaystyle
H^{m}(\mathbb{T}_h)^q := \{ v \in L^2(\Omega)^q: \ v|_{T} \in H^{1}(T)^q \quad \forall T \in \mathbb{T}_h  \}, \quad m \in \mathbb{N}
\end{align*}
with a norm
\begin{align*}
\displaystyle
| v |_{H^{m}(\mathbb{T}_h)^q} &:= \left( \sum_{T \in \mathbb{T}_h} | v |^2_{H^{m}(T)^q } \right)^{\frac{1}{2}}.
\end{align*}
When $q=1$, we denote $H^{m}(\mathbb{T}_h) := H^{m}(\mathbb{T}_h)^1 $. Let $\varphi \in H^{1}(\mathbb{T}_h)$. Suppose that $F \in \mathcal{F}_h^i$ with $F = T_1 \cap T_2$, $T_1,T_2 \in \mathbb{T}_h$. Set $\varphi_1 := \varphi{|_{T_1}}$ and $\varphi_2 := \varphi{|_{T_2}}$. Set two nonnegative real numbers $\omega_{T_1,F}$ and $\omega_{T_2,F}$ such that
\begin{align*}
\displaystyle
\omega_{T_1,F} + \omega_{T_2,F} = 1.
\end{align*}
The jump and the skew-weighted average of $\varphi$ across $F$ is then defined as
\begin{align*}
\displaystyle
[\![\varphi]\!] := [\! [ \varphi ]\!]_F := \varphi_1 - \varphi_2, \quad  \{\! \{ \varphi\} \! \}_{\overline{\omega}} :=  \{\! \{ \varphi\} \! \}_{\overline{\omega},F} := \omega_{T_2,F} \varphi_1 + \omega_{T_1,F} \varphi_2.
\end{align*}
For a boundary face $F \in \mathcal{F}_h^{\partial}$ with $F = \partial T \cap \partial \Omega$, $[\![\varphi ]\!]_F := \varphi|_{T}$ and $\{\! \{ \varphi \} \!\}_{\overline{\omega}} := \varphi |_{T}$. For any $v \in H^{1}(\mathbb{T}_h)^d$, we use the notation
\begin{align*}
\displaystyle
&[\![v \cdot n]\!] := [\![ v \cdot n ]\!]_F := v_1 \cdot n_1 + v_2 \cdot n_2,  \  \{\! \{ v\} \! \}_{\omega} :=  \{\! \{ v \} \! \}_{\omega,F} := \omega_{T_1,F} v_1 + \omega_{T_2,F} v_2,\\
&[\![v]\!] :=  [\![ v]\!]_F := v_1 - v_2,
\end{align*}
for the jump of the normal component of $v$, the weighted average of $v$, and the junp {of} $v$. For $F \in \mathcal{F}_h$, $[\![v \cdot n]\!] := v \cdot n$, $\{\! \{ v\} \! \}_{\omega} := v$ and $[\![v]\!] :=$, where $n$ is the outward normal.

We define a broken gradient operator as follows. For $\varphi \in H^1(\mathbb{T}_h)$, the broken gradient $\nabla_h:H^1(\mathbb{T}_h) \to L^2(\Omega)^{d}$ is defined by
\begin{align*}
\displaystyle
(\nabla_h \varphi)|_T &:= \nabla (\varphi|_T) \quad \forall T \in \mathbb{T}_h.
\end{align*}

\subsection{Penalty parameters and energy norms} \label{pena=norm}
The following trace inequality on anisotropic meshes is significant in this study. Some references can be found for the proof. Here, we follow Ern and Guermond \cite[Lemma 12.15]{ErnGue21a}.

\begin{lem}[Trace inequality] \label{lem=trace}
Let  $T \subset \mathbb{R}^d$ be a simplex. There exists a positive constant $c$ such that for any $\varphi \in H^{1}(T)$, $F \in \mathcal{F}_{T}$, and $h$,
\begin{align}
\displaystyle
\| \varphi \|_{L^2(F)}
\leq c \ell_{T,F}^{- \frac{1}{2}} \left( \| \varphi \|_{L^2(T)} + h_{T}^{\frac{1}{2}}  \| \varphi \|_{L^2(T)}^{\frac{1}{2}} | \varphi |_{H^1(T)}^{\frac{1}{2}} \right),\label{trace}
\end{align}
where $\ell_{T,F} := \frac{d! |T|_d}{|F|_{d-1}}$ denotes the distance of the vertex of $T$ opposite to $F$ to the face. Furthermore, there exists a positive constant $c$ such that for any $v = (v^{(1)}, \ldots,v^{(d)})^T \in H^{1}(T)^d$, $F \in \mathcal{F}_{T}$, and $h$,
\begin{align}
\displaystyle
\| v \|_{L^2(F)^d}
\leq c \ell_{T,F}^{- \frac{1}{2}} \left( \| v \|_{L^2(T)^d} + h_{T}^{\frac{1}{2}}  \| v \|_{L^2(T)^d}^{\frac{1}{2}} | v |_{H^1(T)^d}^{\frac{1}{2}} \right),\label{trace=vec}
\end{align}
\end{lem}

\begin{pf*}
A proof is found in \cite[Lemma 1]{Ish24}.
\qed
\end{pf*}

Deriving an appropriate penalty term is essential in discontinuous Galerkin methods (dG) on anisotropic meshes. The use of weighted averages gives robust dG schemes for various problems; see \cite{DonGeo22,Ish24,PieErn12}. For any $v \in H^{1}(\mathbb{T}_h)^d$ and $\varphi \in H^{1}(\mathbb{T}_h)$,
\begin{align}
\displaystyle
[\![ ((v \varphi) \cdot n ]\!]_F
&=  \{\! \{ v \} \! \}_{\omega,F} \cdot n_F [\! [ \varphi  ]\!]_F + [\![ v \cdot n ]\!]_F \{\! \{ \varphi  \} \! \}_{\overline{\omega},F}. \label{jump=ave}
\end{align}
For example, if $u \in H^1_0(\Omega) \cap W^{2,1}(\Omega)$, setting $v := - \nabla u$, we have $ [\![ v \cdot n ]\!]_F = 0$ for all $F \in \mathcal{F}_h^i$, see \cite[Lemma 4.3]{PieErn12}. Using the trace (Lemma \ref{lem=trace}) and the H\"older inequalities, the weighted average gives the following estimate for the term $ \{\! \{ v \}\! \}_{\omega,F} \cdot n_F [\![ \varphi  ]\!]_F$.
\begin{align}
\displaystyle
&\int_{F} \left|   \{\! \{ v \} \! \}_{\omega,F} \cdot n_F [\! [ \varphi ]\!]_F  \right| ds \nonumber\\
&\leq c h^{\beta}  \left( \| v|_{T_1} \|_{H^1(T_1)^d}^2 + \| v|_{T_2} \|_{H^1(T_2)^d}^2 \right)^{\frac{1}{2}} \nonumber\\
&\quad \times \left( h^{- 2 \beta} \omega_{T_1,F}^2 { \ell_{T_1,F}}^{-1} + h^{- 2 \beta}\omega_{T_2,F}^2 {\ell_{T_2,F}}^{-1} \right)^{\frac{1}{2}} \left \|  [\![ \varphi  ]\!] \right\|_{L^2(F)}, {\label{new=2=9}}
\end{align}
{where $\ell_{T_1,F}$ and $\ell_{T_2,F}$ are defined in the inequality \eqref{trace=vec}. The weights $\omega_{T_1,F}$, $\omega_{T_2,F}$ and $\beta$ are nonnegative real numbers chosen latter on.} A choice for the weighted parameters is such that for $F \in \mathcal{F}_h^i$ with $F = T_{1} \cap T_{2}$, $T_{1},T_{2} \in \mathbb{T}_h$,
\begin{align}
\displaystyle
 \omega_{T_i,F} :=  \frac{\sqrt{\ell_{T_i,F}}}{\sqrt{\ell_{T_1,F}} + \sqrt{\ell_{T_2,F}}}, \quad i=1,2. \label{weight}
\end{align}
Then, the associated penalty parameter is defined as
\begin{align}
\displaystyle
 \frac{2 h^{- 2 \beta}}{(\sqrt{ \ell_{T_1,F}} + \sqrt{\ell_{T_2,F}})^2}. \label{weight_pena}
\end{align}
Let $F \in \mathcal{F}_h^i$ with $F = T_1 \cap T_2$, $T_1,T_2 \in \mathbb{T}_h$ be an interior face and $F \in \mathcal{F}_h^{\partial}$ with $F = \partial T_{\partial} \cap \partial \Omega$, $T_{\partial} \in \mathbb{T}_h$ a boundary face. A new penalty parameter $\kappa_F$ for the WOPSIP method is defined as follows using \eqref{weight_pena} with $\beta =1$. 
\begin{align}
\displaystyle
\kappa_F :=
\begin{cases}
\displaystyle
h^{-2} \left( \sqrt{ \ell_{T_1,F}} + \sqrt{\ell_{T_2,F}} \right)^{-2} \quad \text{if $F \in \mathcal{F}_h^i$},\\
\displaystyle
h^{-2} \ell_{T_{\partial},F}^{-1} \quad \text{if $F \in \mathcal{F}_h^{\partial}$}.
\end{cases} \label{penalty0}
\end{align}
For the RSIP method and the discrete Poincar\'e inequality, we use the following parameter.
\begin{align}
\displaystyle
\kappa_{F*} :=
\begin{cases}
\displaystyle
\left( \sqrt{ \ell_{T_1,F}} + \sqrt{\ell_{T_2,F}} \right)^{-2} \quad \text{if $F \in \mathcal{F}_h^i$},\\
\displaystyle
\ell_{T_{\partial},F}^{-1} \quad \text{if $F \in \mathcal{F}_h^{\partial}$}.
\end{cases} \label{penalty1}
\end{align}

For any $F \in \mathcal{F}_h$, we define the $L^2$-projection $\Pi_F^{0}: L^2(F) \to \mathbb{P}^{0}(F)$ by
\begin{align*}
\displaystyle
\int_F (\Pi_F^{0} \varphi - \varphi)   ds = 0 \quad \forall \varphi \in L^2(F).
\end{align*}
We then define the following norms for any $v \in H^1(\mathbb{T}_h)$.
\begin{align*}
\displaystyle
| v |_{wop} := \left( | v |_{H^1(\mathbb{T}_h)}^2 + | v |_{jwop}^2 \right)^{\frac{1}{2}}
\end{align*}
with the jump seminorm
\begin{align*}
\displaystyle
 | v |_{jwop} := \left( \sum_{F \in \mathcal{F}_h} \kappa_F \| \Pi_F^{0} [\![ v ]\!] \|_{L^2(F)}^2  \right)^{\frac{1}{2}}
\end{align*}
and $\kappa_F$ defined as in \eqref{penalty0};
\begin{align*}
\displaystyle
| v |_{rdg} := \left( | v |_{H^1(\mathbb{T}_h)}^2 + | v |_{jrdg}^2 \right)^{\frac{1}{2}}
\end{align*}
with
\begin{align*}
\displaystyle
 | v |_{jrdg} := \left( \sum_{F \in \mathcal{F}_h} \kappa_{F*} \| \Pi_F^{0} [\![ v ]\!] \|_{L^2(F)}^2  \right)^{\frac{1}{2}}
\end{align*}
and $\kappa_{F*}$ defined as in \eqref{penalty1}. For any $v \in H^1(\mathbb{T}_h)$, $ | v |_{jrdg} \leq  | v |_{jwop}$ for $h \leq 1$. The norm $| \cdot |_{wop}$ is used for analysis of the WOPSIP method, while the norm $| \cdot |_{rdg}$ is used for analysis of the RSIP method and the discrete Poincar\'e inequality (Lemma \ref{Poin=lem6}).

Furthermore, for any $k \in \mathbb{N}_0$, let $\mathbb{P}^k(T)$ and $\mathbb{P}^k(F)$ be spaces of polynomials with degree at most $k$ in $T$ and $F$, respectively. 

\subsection{Edge characterisation on a simplex, a geometric parameter, and a condition} \label{element=cond}
We impose edge characterisation on a simplex to analyse anisotropic error estimates.

\subsubsection{Reference elements} \label{reference}
We now define the reference elements $\widehat{T} \subset \mathbb{R}^d$.

\subsubsection*{Two-dimensional case} \label{reference2d}
Let $\widehat{T} \subset \mathbb{R}^2$ be a reference triangle with vertices $\hat{p}_1 := (0,0)^T$, $\hat{p}_2 := (1,0)^T$, and $\hat{p}_3 := (0,1)^T$. 

\subsubsection*{Three-dimensional case} \label{reference3d}
In the three-dimensional case, we consider the following two cases: (\roman{sone}) and (\roman{stwo}); see Condition \ref{cond2}.

Let $\widehat{T}_1$ and $\widehat{T}_2$ be reference tetrahedra with the following vertices:
\begin{description}
   \item[(\roman{sone})] $\widehat{T}_1$ has vertices $\hat{p}_1 := (0,0,0)^T$, $\hat{p}_2 := (1,0,0)^T$, $\hat{p}_3 := (0,1,0)^T$, and $\hat{p}_4 := (0,0,1)^T$;
 \item[(\roman{stwo})] $\widehat{T}_2$ has vertices $\hat{p}_1 := (0,0,0)^T$, $\hat{p}_2 := (1,0,0)^T$, $\hat{p}_3 := (1,1,0)^T$, and $\hat{p}_4 := (0,0,1)^T$.
\end{description}
Therefore, we set $\widehat{T} \in \{ \widehat{T}_1 , \widehat{T}_2 \}$. 
%We only use case  (\roman{sone}) to show the RT interpolation error estimates, see Theorem \ref{thr3}. 
Note that the case (\roman{sone}) is called \textit{the regular vertex property}, see \cite{AcoDur99}.

\subsubsection{Affine mappings} \label{Affinedef}
We introduced a new strategy proposed in \cite[Section 2]{IshKobTsu21c} to use anisotropic mesh partitions. We construct two affine simplex $\widetilde{T} \subset \mathbb{R}^d$, we construct two affine mappings $\Phi_{\widetilde{T}}: \widehat{T} \to \widetilde{T}$ and $\Phi_{T}: \widetilde{T} \to T$. First, we define the affine mapping $\Phi_{\widetilde{T}}: \widehat{T} \to \widetilde{T}$ as
\begin{align}
\displaystyle
\Phi_{\widetilde{T}}: \widehat{T} \ni \hat{x} \mapsto \tilde{x} := \Phi_{\widetilde{T}}(\hat{x}) := {A}_{\widetilde{T}} \hat{x} \in  \widetilde{T}, \label{aff=1}
\end{align}
where ${A}_{\widetilde{T}} \in \mathbb{R}^{d \times d}$ is an invertible matrix. Section \ref{sec221} provides the details. We then define the affine mapping $\Phi_{T}: \widetilde{T} \to T$  as follows:
\begin{align}
\displaystyle
\Phi_{T}: \widetilde{T} \ni \tilde{x} \mapsto x := \Phi_{T}(\tilde{x}) := {A}_{T} \tilde{x} + b_{T} \in T, \label{aff=2}
\end{align}
where $b_{T} \in \mathbb{R}^d$ is a vector and ${A}_{T} \in O(d)$ is the rotation and mirror imaging matrix. Section \ref{sec322} provides the details. We define the affine mapping $\Phi: \widehat{T} \to T$ as
\begin{align*}
\displaystyle
\Phi := {\Phi}_{T} \circ {\Phi}_{\widetilde{T}}: \widehat{T} \ni \hat{x} \mapsto x := \Phi (\hat{x}) =  ({\Phi}_{T} \circ {\Phi}_{\widetilde{T}})(\hat{x}) = {A} \hat{x} + b_{T} \in T, 
\end{align*}
where ${A} := {A}_{T} {A}_{\widetilde{T}} \in \mathbb{R}^{d \times d}$.

\subsubsection{Construct mapping $\Phi_{\widetilde{T}}: \widehat{T} \to \widetilde{T}$} \label{sec221} 
We consider affine mapping \eqref{aff=1}. We define the matrix $ {A}_{\widetilde{T}} \in \mathbb{R}^{d \times d}$ as follows: We first define the diagonal matrix as
\begin{align}
\displaystyle
\widehat{A} :=  \diag (h_1,\ldots,h_d), \quad h_i \in \mathbb{R}_+ \quad \forall i,\label{aff=3}
\end{align}
where $\mathbb{R}_+$ denotes the set of positive real numbers.

For $d=2$, we define the regular matrix $\widetilde{A} \in \mathbb{R}^{2 \times 2}$ as:
\begin{align}
\displaystyle
\widetilde{A} :=
\begin{pmatrix}
1 & s \\
0 & t \\
\end{pmatrix}, \label{aff=4}
\end{align}
with parameters
\begin{align*}
\displaystyle
s^2 + t^2 = 1, \quad t \> 0.
\end{align*}
For reference element $\widehat{T}$, let $\mathfrak{T}^{(2)}$ be a family of triangles.
\begin{align*}
\displaystyle
\widetilde{T} &= \Phi_{\widetilde{T}}(\widehat{T}) = {A}_{\widetilde{T}} (\widehat{T}), \quad {A}_{\widetilde{T}} := \widetilde {A} \widehat{A}
\end{align*}
with vertices $\tilde{p}_1 := (0,0)^T$, $\tilde{p}_2 := (h_1,0)^T$, and $\tilde{p}_3 :=(h_2 s , h_2 t)^T$. Then,  $h_1 = |\tilde{p}_1 - \tilde{p}_2| \> 0$ and $h_2 = |\tilde{p}_1 - \tilde{p}_3| \> 0$. 

For $d=3$, we define the regular matrices $\widetilde{A}_1, \widetilde{A}_2 \in \mathbb{R}^{3 \times 3}$ as 
\begin{align}
\displaystyle
\widetilde{A}_1 :=
\begin{pmatrix}
1 & s_1 & s_{21} \\
0 & t_1  & s_{22}\\
0 & 0  & t_2\\
\end{pmatrix}, \
\widetilde{A}_2 :=
\begin{pmatrix}
1 & - s_1 & s_{21} \\
0 & t_1  & s_{22}\\
0 & 0  & t_2\\
\end{pmatrix} \label{aff=5}
\end{align}
with parameters
\begin{align*}
\displaystyle
\begin{cases}
s_1^2 + t_1^2 = 1, \ s_1 \> 0, \ t_1 \> 0, \ h_2 s_1 \leq h_1 / 2, \\
s_{21}^2 + s_{22}^2 + t_2^2 = 1, \ t_2 \> 0, \ h_3 s_{21} \leq h_1 / 2.
\end{cases}
\end{align*}
Therefore, we set $\widetilde{A} \in \{ \widetilde{A}_1 , \widetilde{A}_2 \}$. For the reference elements $\widehat{T}_i$, $i=1,2$, let $\mathfrak{T}_i^{(3)}$ and $i=1,2$ be a family of tetrahedra.
\begin{align*}
\displaystyle
\widetilde{T}_i &= \Phi_{\widetilde{T}_i} (\widehat{T}_i) =  {A}_{\widetilde{T}_i} (\widehat{T}_i), \quad {A}_{\widetilde{T}_i} := \widetilde {A}_i \widehat{A}, \quad i=1,2,
\end{align*}
with vertices
\begin{align*}
\displaystyle
&\tilde{p}_1 := (0,0,0)^T, \ \tilde{p}_2 := (h_1,0,0)^T, \ \tilde{p}_4 := (h_3 s_{21}, h_3 s_{22}, h_3 t_2)^T, \\
&\begin{cases}
\tilde{p}_3 := (h_2 s_1 , h_2 t_1 , 0)^T \quad \text{for case (\roman{sone})}, \\
\tilde{p}_3 := (h_1 - h_2 s_1, h_2 t_1,0)^T \quad \text{for case (\roman{stwo})}.
\end{cases}
\end{align*}
Subsequently, $h_1 = |\tilde{p}_1 - \tilde{p}_2| \> 0$, $h_3 = |\tilde{p}_1 - \tilde{p}_4| \> 0$, and
\begin{align*}
\displaystyle
h_2 =
\begin{cases}
|\tilde{p}_1 - \tilde{p}_3| \> 0  \quad \text{for case (\roman{sone})}, \\
|\tilde{p}_2 - \tilde{p}_3| \> 0  \quad \text{for case (\roman{stwo})}.
\end{cases}
\end{align*}

\subsubsection{Construct mapping $\Phi_{T}: \widetilde{T} \to T$}  \label{sec322}
We determine the affine mapping \eqref{aff=2} as follows: Let ${T} \in \mathbb{T}_h$ have vertices ${p}_i$ ($i=1,\ldots,d+1$). Let $b_{T} \in \mathbb{R}^d$ be the vector and ${A}_{T} \in O(d)$ be the rotation and mirror imaging matrix such that
\begin{align*}
\displaystyle
p_{i} = \Phi_T (\tilde{p}_i) = {A}_{T} \tilde{p}_i + b_T, \quad i \in \{1, \ldots,d+1 \},
\end{align*}
where the vertices $p_{i}$ ($i=1,\ldots,d+1$) satisfy the following conditions:

\begin{Cond}[Case in which $d=2$] \label{cond1}
Let ${T} \in \mathbb{T}_h$ with the vertices ${p}_i$ ($i=1,\ldots,3$). We assume that $\overline{{p}_2 {p}_3}$ is the longest edge of ${T}$; i.e., $ h_{{T}} := |{p}_2 - {p}_ 3|$. We set  $h_1 = |{p}_1 - {p}_2|$ and $h_2 = |{p}_1 - {p}_3|$. We then assume that $h_2 \leq h_1$. Note that ${h_1 \approx h_T}$. %$h_1 = \mathcal{O}(h_{{T}})$. 
\end{Cond}

\begin{Cond}[Case in which $d=3$] \label{cond2}
Let ${T} \in \mathbb{T}_h$ with the vertices ${p}_i$ ($i=1,\ldots,4$). Let ${L}_i$ ($1 \leq i \leq 6$) be the edges of ${T}$. We denote by ${L}_{\min}$  the edge of ${T}$ with the minimum length; i.e., $|{L}_{\min}| = \min_{1 \leq i \leq 6} |{L}_i|$. We set $h_2 := |{L}_{\min}|$ and assume that 
\begin{align*}
\displaystyle
&\text{the endpoints of ${L}_{\min}$ are either $\{ {p}_1 , {p}_3\}$ or $\{ {p}_2 , {p}_3\}$}.
\end{align*}
Among the four edges that share an endpoint with ${L}_{\min}$, we take the longest edge ${L}^{({\min})}_{\max}$. Let ${p}_1$ and ${p}_2$ be the endpoints of edge ${L}^{({\min})}_{\max}$. We thus have that
\begin{align*}
\displaystyle
h_1 = |{L}^{(\min)}_{\max}| = |{p}_1 - {p}_2|.
\end{align*}
We consider cutting $\mathbb{R}^3$ with the plane that contains the midpoint of edge ${L}^{(\min)}_{\max}$ and is perpendicular to the vector ${p}_1 - {p}_2$. Thus, we have two cases: 
\begin{description}
  \item[(Type \roman{sone})] ${p}_3$ and ${p}_4$  belong to the same half-space;
  \item[(Type \roman{stwo})] ${p}_3$ and ${p}_4$  belong to different half-spaces.
\end{description}
In each case, we set
\begin{description}
  \item[(Type \roman{sone})] ${p}_1$ and ${p}_3$ as the endpoints of ${L}_{\min}$, that is, $h_2 =  |{p}_1 - {p}_3| $;
  \item[(Type \roman{stwo})] ${p}_2$ and ${p}_3$ as the endpoints of ${L}_{\min}$, that is, $h_2 =  |{p}_2 - {p}_3| $.
\end{description}
Finally, we set $h_3 = |{p}_1 - {p}_4|$. Note that we implicitly assume that ${p}_1$ and ${p}_4$ belong to the same half-space. In addition, note that ${h_1 \approx h_T}$. %$h_1 = \mathcal{O}(h_{{T}})$.
\end{Cond}

\subsection{Additional notation and assumption} \label{addinot}
We define vectors ${r}_n \in \mathbb{R}^d$, $n=1,\ldots,d$ as follows. If $d=2$,
\begin{align*}
\displaystyle
{r}_1 := \frac{p_2 - p_1}{|p_2 - p_1|}, \quad {r}_2 := \frac{p_3 - p_1}{|p_3 - p_1|},
\end{align*}
and if $d=3$,
\begin{align*}
\displaystyle
&{r}_1 := \frac{p_2 - p_1}{|p_2 - p_1|}, \quad {r}_3 := \frac{p_4 - p_1}{|p_4 - p_1|}, \quad
\begin{cases}
\displaystyle
{r}_2 := \frac{p_3 - p_1}{|p_3 - p_1|}, \quad \text{for case (\roman{sone})}, \\
\displaystyle
{r}_2 := \frac{p_3 - p_2}{|p_3 - p_2|} \quad \text{for case (\roman{stwo})}.
\end{cases}
\end{align*}
For a sufficiently smooth function $\varphi$ and vector function $v := (v_{1},\ldots,v_{d})^T$, we define the directional derivative as, for $i \in \{ 1, \ldots,d \}$,
\begin{align*}
\displaystyle
\frac{\partial \varphi}{\partial {r_i}} &:= ( {r}_i \cdot  \nabla_{x} ) \varphi = \sum_{i_0=1}^d ({r}_i)_{i_0} \frac{\partial \varphi}{\partial x_{i_0}^{}}, \\
\frac{\partial v}{\partial r_i} &:= \left(\frac{\partial v_{1}}{\partial r_i}, \ldots, \frac{\partial v_{d}}{\partial r_i} \right)^T 
= ( ({r}_i  \cdot \nabla_{x}) v_{1}, \ldots, ({r}_i  \cdot \nabla_{x} ) v_{d} )^T.
\end{align*}
For a multi-index $\beta = (\beta_1,\ldots,\beta_d) \in \mathbb{N}_0^d$, we use the notation
\begin{align*}
\displaystyle
\partial^{\beta} \varphi := \frac{\partial^{|\beta|} \varphi}{\partial x_1^{\beta_1} \ldots \partial x_d^{\beta_d}}, \quad \partial^{\beta}_{r} \varphi := \frac{\partial^{|\beta|} \varphi}{\partial r_1^{\beta_1} \ldots \partial r_d^{\beta_d}}, \quad h^{\beta} :=  h_{1}^{\beta_1} \cdots h_{d}^{\beta_d}.
\end{align*}
Note that $\partial^{\beta} \varphi \neq  \partial^{\beta}_{r} \varphi$.

We proposed a geometric parameter $H_{T}$ in a prior work \cite{IshKobTsu21a}.
 \begin{defi} \label{defi1}
 The parameter $H_{{T}}$ is defined as
\begin{align*}
\displaystyle
H_{{T}} := \frac{\prod_{i=1}^d h_i}{|{T}|_d} h_{{T}}.
\end{align*}
\end{defi}
We introduce the geometric condition proposed in \cite{IshKobTsu21a}, which is equivalent to the maximum-angle condition \cite{IshKobSuzTsu21d}.

\begin{assume} \label{neogeo=assume}
A family of meshes $\{ \mathbb{T}_h\}$ has a semi-regular property if there exists $\gamma_0 \> 0$ such that
\begin{align}
\displaystyle
\frac{H_{T}}{h_{T}} \leq \gamma_0 \quad \forall \mathbb{T}_h \in \{ \mathbb{T}_h \}, \quad \forall T \in \mathbb{T}_h. \label{NewGeo}
\end{align}
\end{assume}

\subsection{Piola transformations}
The Piola transformation $\Psi : L^1(\widehat{T})^d \to L^1({T})^d$ is defined as
\begin{align*}
\displaystyle
\Psi :  L^1(\widehat{T})^d  &\to  L^1({T})^d \\
\hat{v} &\mapsto v(x) :=  \Psi(\hat{v})(x) = \frac{1}{\det(A)} A \hat{v}(\hat{x}).
\end{align*}

\subsection{Finite element spaces and anisotropic interpolation error estimates}

\subsubsection{Finite element spaces}
For $s \in \mathbb{N}_0$, we define a discontinuous finite element space as
\begin{align*}
\displaystyle
P_{dc,h}^{s} &:= \left\{ p_h \in L^2(\Omega); \ p_h|_{T} \circ {\Phi} \in \mathbb{P}^{s}(\widehat{T}) \quad \forall T \in \mathbb{T}_h  \right\}.
\end{align*}
Let $Ne$ be the number of elements included in the mesh $\mathbb{T}_h$. Thus, we write $\mathbb{T}_h = \{ T_j\}_{j=1}^{Ne}$.

Let the points $\{ P_{T_j,1}, \ldots, P_{T_j,d+1} \}$ be the vertices of the simplex $T_j \in \mathbb{T}_h$ for $j \in \{1, \ldots , Ne \}$. Let $F_{T_j,i}$ be the face of $T_j$ opposite $P_{T_j,i}$ for $i \in \{ 1, \ldots , d+1\}$. We set $P := \mathbb{P}^1$, and take a set ${\Sigma}_{T_j} := \{ {\chi}^{CR}_{T_j,i} \}_{1 \leq i \leq d+1}$ of linear forms with its components such that for any $p \in \mathbb{P}^1$.
\begin{align}
\displaystyle
{\chi}^{CR}_{T_j,i}({p}) := \frac{1}{| {F}_{T_j,i} |_{d-1}} \int_{{F}_{T_j,i}} {p} d{s} \quad \forall i \in \{ 1, \ldots,d+1 \}. \label{CR1}
\end{align}
For each $j \in \{1, \ldots ,Ne \}$, the triple $\{ T_j ,  \mathbb{P}^1 , \Sigma_{T_j} \}$ is a finite element. Using the barycentric coordinates $ \{{\lambda}_{T_j,i} \}_{i=1}^{d+1}: \mathbb{R}^d \to \mathbb{R}$ on the reference element, the nodal basis functions associated with the degrees of freedom by \eqref{CR1} are defined as
\begin{align}
\displaystyle
{\theta}^{CR}_{T_j,i}({x}) := d \left( \frac{1}{d} - {\lambda}_{T_j,i} ({x}) \right) \quad \forall i \in \{ 1, \ldots ,d+1 \}. \label{CR2}
\end{align}
For $j \in \{1, \ldots ,Ne \}$ and $i \in \{ 1, \ldots , d+1\}$, we define the function $\phi_{j(i)}$ as
\begin{align}
\displaystyle
\phi_{j(i)}(x) :=
\begin{cases}
\theta_{T_j,i}^{CR}(x), \quad \text{$x \in T_j$}, \\
0, \quad \text{$x \notin T_j$}.
\end{cases} \label{CR5}
\end{align}
We define a discontinuous finite element space as
\begin{align}
\displaystyle
V_{dc,h}^{CR} &:= \left\{ \sum_{j=1}^{Ne} \sum_{i=1}^{d+1} c_{j(i)} \phi_{j(i)}; \  c_{j(i)} \in \mathbb{R}, \ \forall i,j \right\} \subset P_{dc,h}^1 \label{CR6}
\end{align}
with a norm $|\varphi_h|_{V_{dc,h}^{CR}} := |\varphi_h|_{h}$ for any $\varphi_h \in V_{dc,h}^{CR}$. 
\begin{comment}
We further define a discontinuous space as
\begin{align}
\displaystyle
V_{dc,h,0}^{CR} &:= \left\{ \varphi_h \in V_{dc,h}^{CR}: \ \int_{F} \varphi_h ds = 0 \ \forall F \in \mathcal{F}_h^{\partial} \right\}.
\end{align}
\end{comment}

For $T_j \in \mathbb{T}_h$, $j \in \{ 1, \ldots ,Ne \}$, we define the local RT polynomial space as follows.
\begin{align}
\displaystyle
\mathbb{RT}^0(T_j) := \mathbb{P}^0(T_j)^d + x \mathbb{P}^0(T_j), \quad x \in \mathbb{R}^d. \label{RT1}
\end{align}
For $p \in \mathbb{RT}^0(T_j)$, the local degrees of freedom are defined as
\begin{align}
\displaystyle
{\chi}^{RT}_{T_j,i}({p}) := \int_{{F}_{T_j,i}} {p} \cdot n_{T_j,i} d{s} \quad \forall i \in \{ 1, \ldots ,d+1 \}, \label{RT2}
\end{align}
where $n_{T_j,i}$ is a fixed unit normal to ${F}_{T_j,i}$. Setting $\Sigma_{T_j}^{RT} := \{ {\chi}^{RT}_{T_j,i} \}_{1 \leq i \leq d+1}$, the triple $\{ T_j ,  \mathbb{RT}^0 , \Sigma_{T_j}^{RT} \}$ a finite element. The local shape functions are 
\begin{align}
\displaystyle
\theta_{T_j,i}^{RT}(x) := \frac{\iota_{{F_{T_j,i}},{T}_j}}{d |T_j|_d} (x - P_{T_j,i}) \quad \forall i \in \{ 1, \ldots , d+1 \}, \label{RT3}
\end{align}
where $\iota_{{F_{T_j,i}},{T}_j} := 1$ if ${n}_{T_{j,i}}$ points outwards, and $ - 1$ otherwise \cite{ErnGue21a}. We define a discontinuous RT finite element space as follows.
\begin{align}
\displaystyle
V_{dc,h}^{RT} &:= \{ v_h \in L^1(\Omega)^d : \  \Psi^{-1}(v_h|_{T_j}) \in \mathbb{RT}^0(\widehat{T}) \quad \forall T_j \in \mathbb{T}_h, \  j \in \{ 1, \ldots , Ne\} \}. \label{RT7}
\end{align}

\subsubsection{Discontinuous space and the $L^2$-orthogonal projection}
For $T_j \in \mathbb{T}_h$, $j \in \{ 1, \ldots , Ne \}$, let $\Pi_{T_j}^0 : L^2(T_j) \to \mathbb{P}^0$ be the $L^2$-orthogonal projection defined as
\begin{align*}
\displaystyle
\Pi_{T_j}^0 \varphi := \frac{1}{|T_j|_d} \int_{{T_j}} \varphi dx \quad \forall \varphi \in L^2(T_j).
\end{align*}
The following theorem gives an anisotropic error estimate of the projection $\Pi_{T_j}^0$. 
\begin{thr} \label{thr1}
For any $\hat{\varphi} \in H^{1}(\widehat{T})$ with ${\varphi} := \hat{\varphi} \circ {\Phi}^{-1}$,
\begin{align}
\displaystyle
\| \Pi_{T_j}^0 \varphi - \varphi \|_{L^2(T_j)} \leq c \sum_{i=1}^d h_i \left\| \frac{\partial \varphi}{\partial r_i} \right\|_{L^{2}(T_j)}. \label{L2ortho}
\end{align}
\end{thr}

\begin{pf*}
A proof is found in \cite[Theorem 2]{Ish24} and \cite[Theorem 2]{Ish24b}.
\qed
\end{pf*}

We also define the global interpolation $\Pi_h^0$ to the space $P_{dc,h}^{0}$ as
\begin{align*}
\displaystyle
(\Pi_h^0 \varphi)|_{T_j} := \Pi_{{T_j}}^0 (\varphi|_{T_j}) \ \forall T_j \in \mathbb{T}_h, \ j \in \{ 1, \ldots ,Ne\}, \ \forall \varphi \in L^2(\Omega).
\end{align*}

\subsubsection{Discontinuous CR finite element interpolation operator}
For $T_j \in \mathbb{T}_h$, $j \in \{ 1, \ldots , Ne \}$, let $I_{T_j}^{CR} : H^{1}(T_j) \to \mathbb{P}^1(T_j)$ be the CR interpolation operator such that for any $\varphi \in H^{1}(T_j)$, 
\begin{align*}
\displaystyle
I_{T_j}^{CR}: H^{1}(T_j) \ni \varphi \mapsto I_{T_j}^{CR} \varphi := \sum_{i=1}^{d+1} \left(  \frac{1}{| {F}_{T_j,i} |_{d-1}} \int_{{F}_{T_j,i}} {\varphi} d{s} \right) \theta_{T_j,i}^{CR} \in \mathbb{P}^1(T_j). 
\end{align*}
We then present estimates of the anisotropic CR interpolation error. 

\begin{thr} \label{DGCE=thr3}
For $j \in \{ 1, \ldots , Ne\}$,
\begin{align}
\displaystyle
|I_{T_j}^{CR} \varphi - \varphi |_{H^{1}({T}_j)} &\leq c \sum_{i=1}^d {h}_i \left\| \frac{\partial }{\partial r_i} \nabla \varphi \right \|_{L^2(T_j)^d} \quad \forall {\varphi} \in H^{2}({T}_j), \label{CR4} \\
\|I_{T_j}^{CR} \varphi - \varphi \|_{L^2(T_j)}
&\leq c  \sum_{|\varepsilon| = 2} h^{\varepsilon} \left\| \partial_{r}^{\varepsilon} \varphi  \right\|_{L^{2}(T)} \quad \forall {\varphi} \in H^{2}({T}_j). \label{L2=CR4}
\end{align}
\end{thr}

\begin{pf*}
The proof of \eqref{CR4} is found in \cite[Theorem 3]{Ish24} and \cite[Theorem 3]{Ish24b}.

Let ${\varphi} \in H^{2}({T}_j)$ for $j \in \{ 1, \ldots , Ne\}$.  Using the scaling argument yields
\begin{align}
\displaystyle
\| I_{T_j}^{CR} \varphi - \varphi \|_{L^2(T)}
\leq c |\det({A})|^{\frac{1}{2}} \| I_{\widehat{T}}^{CR} \hat{\varphi} - \hat{\varphi} \|_{L^2(\widehat{T})}. \label{th2=1}
\end{align}
For any $\hat{\eta} \in \mathbb{P}^{1}$, we have that 
\begin{align}
\displaystyle
\| I_{\widehat{T}}^{CR} \hat{\varphi} - \hat{\varphi} \|_{L^2(\widehat{T})}
&\leq \| I_{\widehat{T}}^{CR}  (\hat{\varphi} - \hat{\eta}) \|_{L^2(\widehat{T})}  + \| \hat{\eta} - \hat{\varphi} \|_{L^2(\widehat{T})},  \label{th2=2}
\end{align}
because $I_{\widehat{T}}^{CR} \hat{\eta} = \hat{\eta}$. Using the trace inequality on the reference element, 
\begin{align}
\displaystyle
\| I_{\widehat{T}}^{CR}  (\hat{\varphi} - \hat{\eta}) \|_{L^2(\widehat{T})} 
&\leq c \| \hat{\varphi} - \hat{\eta} \|_{H^1(\widehat{T})}. \label{th2=3}
\end{align}
Based on \eqref{th2=1}, \eqref{th2=2} and \eqref{th2=3}, we have that
\begin{align}
\displaystyle
\| I_{T_j}^{CR} \varphi - \varphi \|_{L^2(T)}
\leq c |\det({A})|^{\frac{1}{2}} \inf_{\hat{\eta} \in \mathbb{P}^{1}} \| \hat{\varphi} - \hat{\eta} \|_{H^1(\widehat{T})}. \label{th2=4}
\end{align}
From the Bramble--Hilbert lemma (refer to \cite[Lemma 4.3.8]{BreSco08}),  $\hat{\eta}_{\beta} \in \mathbb{P}^{1}$ exists such that for any $\hat{\varphi} \in H^{2}(\widehat{T})$,
\begin{align}
\displaystyle
| \hat{\varphi} - \hat{\eta}_{\beta} |_{H^{s}(\widehat{T})} \leq C^{BH}(\widehat{T}) |\hat{\varphi}|_{H^{2}(\widehat{T})}, \quad s=0,1. \label{th2=5}
\end{align}
Using the inequality in \cite[Lemma 6]{IshKobTsu21c} with $m=0$, we can estimate inequality \eqref{th2=5} as
\begin{align}
\displaystyle
|\hat{\varphi}|_{H^{2}(\widehat{T})}
\leq c |\det({{A}})|^{-\frac{1}{2}} \sum_{|\varepsilon| = 2} h^{\varepsilon} \left\| \partial_{r}^{\varepsilon} \varphi  \right\|_{L^{2}(T)}.  \label{th2=6}
\end{align}
Using \eqref{th2=4}, \eqref{th2=5} and \eqref{th2=6}, we can deduce the target inequality \eqref{L2=CR4}.
\qed
\end{pf*}

We define a global interpolation operator ${I}_{h}^{CR}: H^{1}(\Omega) \to V_{dc,h}^{CR}$ as
\begin{align}
\displaystyle
({I}_{h}^{CR} \varphi )|_{T_j} = {I}_{T_j}^{CR} (\varphi |_{T_j}), \quad j \in \{ 1, \ldots , Ne\}, \quad \forall \varphi \in H^{1}(\Omega).  \label{CR8}
\end{align}	

\subsubsection{Discontinuous RT finite element interpolation operator} \label{RTsp}
For $T_j \in \mathbb{T}_h$, $j \in \{ 1, \ldots , Ne \}$, let $\mathcal{I}_{T_j}^{RT}: H^{1}(T_j)^d \to \mathbb{RT}^0(T_j)$ be the RT interpolation operator such that for any $v \in H^{1}(T_j)^d$,
\begin{align}
\displaystyle
\mathcal{I}_{T_j}^{RT}: H^{1}(T_j)^d \ni v \mapsto \mathcal{I}_{T_j}^{RT} v := \sum_{i=1}^{d+1} \left(  \int_{{F}_{T_j,i}} {v} \cdot n_{T_j,i} d{s} \right) \theta_{T_j,i}^{RT} \in \mathbb{RT}^0(T_j). \label{RT4}
\end{align}
The following two theorems are divided into the element of (Type \roman{sone}) or the element of (Type \roman{stwo}) in Section \ref{element=cond} when $d=3$.

\begin{thr} \label{DGRT=thr3}
Let $T_j$ be the element with Conditions \ref{cond1} or \ref{cond2} and satisfy (Type \roman{sone}) in Section \ref{element=cond} when $d=3$. For any $\hat{v} \in H^{1}(\widehat{T})^d$ with ${v} = ({v}_1,\ldots,{v}_d)^T := {\Psi} \hat{v}$ and $j \in \{ 1, \ldots ,Ne\}$,
\begin{align}
\displaystyle
\| \mathcal{I}_{T_j}^{RT} v - v \|_{L^2(T_j)^d} 
&\leq  c \left( \frac{H_{T_j}}{h_{T_j}} \sum_{i=1}^d h_i \left \|  \frac{\partial v}{\partial r_i} \right \|_{L^2(T_j)^d} +  h_{T_j} \| \div {v} \|_{L^{2}({T}_j)} \right). \label{RT5}
\end{align}
\end{thr}

\begin{pf*}
A proof is provided in \cite[Theorem 2]{Ish21}.
\qed
\end{pf*}

\begin{thr} \label{DGRT=thr4}
Let $d=3$. Let $T_j$ be an element with Condition \ref{cond2} that satisfies (Type \roman{stwo}) in Section \ref{element=cond}. For $\hat{v} \in H^{1}(\widehat{T})^3$ with ${v} = ({v}_1,v_2,{v}_3)^T := {\Psi} \hat{v}$ and $j \in \{ 1, \ldots ,Ne\}$,
\begin{align}
\displaystyle
&\| \mathcal{I}_{T_j}^{RT} v - v \|_{L^2(T_j)^3} 
\leq c \frac{H_{T_j}}{h_{T_j}} \Biggl(  h_{T_j} |v|_{H^1(T_j)^3} \Biggr). \label{RT6}
\end{align}
\end{thr}

\begin{pf*}
A proof is provided in \cite[Theorem 3]{Ish21}.
\qed
\end{pf*}

We define a global  interpolation operator $\mathcal{I}_{h}^{RT}: H^{1}(\Omega)^d  \to V_{dc,h}^{RT}$ by
\begin{align}
\displaystyle
(\mathcal{I}_{h}^{RT} v )|_{T_j} = \mathcal{I}_{T_j}^{RT} (v |_{T_j}), \quad j \in \{ 1, \ldots , Ne\}, \quad \forall v \in H^{1}(\Omega)^d.  \label{RT8}
\end{align}

\subsection{Relation between the RT interpolation and the $L^2$-projection}
Between the RT interpolation $\mathcal{I}_{h}^{RT}$ and the $L^2$-projection $\Pi_h^0$, the following relation holds.

\begin{lem} \label{lem2}
For $j \in \{ 1, \ldots , Ne\}$,
\begin{align}
\displaystyle
\div (\mathcal{I}_{T_j}^{RT} v) = \Pi_{T_j}^0 (\div v) \quad \forall v \in H^{1}(T_j)^d. \label{RT9a}
\end{align}
By combining \eqref{RT9a}, for any $v \in H^1(\Omega)^d$
\begin{align}
\displaystyle
\div (\mathcal{I}_{h}^{RT} v) = \Pi_h^0 (\div v). \label{RT9b}
\end{align}
\end{lem}

\begin{pf*}
A proof is provided in \cite[Lemma 16.2]{ErnGue21a}.
\qed
\end{pf*}

\subsection{Existing results}

\subsubsection{Important tools}
The following relation plays an important role in the discontinuous Galerkin finite element analysis on anisotropic meshes.

\begin{lem} \label{useful=lem3}
For any $w \in H^1(\Omega)^d$ and $\psi_h \in P_{dc,h}^{1}$,
\begin{align}
\displaystyle
&\int_{\Omega} \left( \mathcal{I}_h^{RT} w \cdot \nabla_h \psi_{h} + \div \mathcal{I}_h^{RT} w  \psi_{h} \right) dx\notag\\
&\quad =  \sum_{F \in \mathcal{F}_h^i} \int_{F}  \{ \! \{ w \} \!\}_{\omega,F} \cdot n_F \Pi_F^0 [\![ \psi_{h} ]\!]_F ds + \sum_{F \in \mathcal{F}_h^{\partial}} \int_{F} ( w \cdot n_F) \Pi_F^0 \psi_{h} ds.  \label{wop=3}
\end{align}
\end{lem}

\begin{pf*}
A proof is provided in \cite[Lemma 3]{Ish24}.
\qed
\end{pf*}

The right-hand terms in \eqref{wop=3} are estimated as follows.

\begin{lem} \label{useful=lem4}
For any $w \in H^1(\Omega)^d$ and $\psi_h \in P_{dc,h}^{1}$,
\begin{align}
\displaystyle
&\left| \sum_{F \in \mathcal{F}_h^i} \int_{F}  \{ \! \{ w \} \!\}_{\omega,F} \cdot n_F \Pi_F^0 [\![ \psi_{h} ]\!]_F ds \right|  \notag\\
&\quad \leq c |\psi_{h}|_{jwop} \left(  h \| w \|_{L^2(\Omega)^d} + h^{\frac{3}{2}} \| w \|_{L^2(\Omega)^d}^{\frac{1}{2}} | w |_{H^1(\Omega)^d}^{\frac{1}{2}} \right), \label{wop=5} \\
&\left| \sum_{F \in \mathcal{F}_h^{\partial}} \int_{F} ( w \cdot n_F) \Pi_F^0 \psi_{h} ds \right| \notag \\
&\quad \leq c |\psi_{h}|_{jwop} \left(  h \| w \|_{L^2(\Omega)^d} + h^{\frac{3}{2}} \| w \|_{L^2(\Omega)^d}^{\frac{1}{2}} | w |_{H^1(\Omega)^d}^{\frac{1}{2}} \right). \label{wop=6}
\end{align}
%where $\mathbb{T}_F$ denotes the set of the simplices in $\mathbb{T}_h$ that share $F$ as a common face.
\end{lem}

\begin{pf*}
A proof is provided in \cite[Lemma 4]{Ish24}.
\qed
\end{pf*}

\begin{lem} \label{useful=lem5}
Let $h \leq 1$. Thus, for any $w \in H^1(\Omega)^d$ and $\psi_h \in P_{dc,h}^{1}$,
\begin{align}
\displaystyle
\left| \sum_{F \in \mathcal{F}_h^i} \int_{F}  \{ \! \{ w \} \!\}_{\omega,F} \cdot n_F \Pi_F^0 [\![ \psi_{h} ]\!]_F ds \right|
&\leq c |\psi_{h}|_{rdg} \| w \|_{H^1(\Omega)^d}, \label{wop=7} \\
\left| \sum_{F \in \mathcal{F}_h^{\partial}} \int_{F} ( w \cdot n_F) \Pi_F^0 \psi_{h} ds \right|
&\leq  c |\psi_{h}|_{rdg}  \| w \|_{H^1(\Omega)^d}. \label{wop=8}
\end{align}
\end{lem}

\begin{pf*}
A proof is provided in \cite[Lemma 4]{Ish24}.
\qed
\end{pf*}

\subsubsection{Discrete Poincar\'e inequality}
The following lemma provides a discrete Poincar\'e inequality. For simplicity, we assume that $\Omega$ is convex.

\begin{lem}[Discrete Poincar\'e inequality] \label{Poin=lem6}
Assume that $\Omega$ is convex. Let $\{ \mathbb{T}_h\}$ be a family of meshes with the semi-regular property (Assumption \ref{neogeo=assume}) and $h \leq 1$. Then, there exists a positive constant $C_{dc}^{P}$ independent of $h$ {but dependent on the maximum angle} such that
\begin{align}
\displaystyle
 \| \psi_{h} \|_{L^2(\Omega)} \leq C_{dc}^{P} | \psi_{h} |_{rdg} \quad \forall \psi_h \in P_{dc,h}^{1}. \label{disPoincare}
\end{align}
\end{lem}

\begin{pf*}
A proof is provided in \cite[Lemma 6]{Ish24}.
\qed
\end{pf*}

\begin{rem} \label{poin=rem1}
For any $f \in L^2(\Omega)$, we set $\ell_h: P_{dc,h}^{1} \to \mathbb{R}$ such that
\begin{align}
\displaystyle
\ell_h(\psi_{h}) := \int_{\Omega} f \psi_{h} dx \quad \forall \psi_{h} \in P_{dc,h}^{1}. \label{lh=rem1}
\end{align}
The H\"older's and the discrete Poincar\'e inequalities yield
\begin{align*}
\displaystyle
|\ell_h(\psi_{h})| \leq C_{dc}^{P} \| f \|_{L^2(\Omega)} | \psi_{h} |_{rdg} \quad \forall \psi_{h} \in P_{dc,h}^{1},
\end{align*}
if $\Omega$ is convex. We are interested in case that $\Omega$ is not convex to prove stability estimates of schemes. In \cite{Bre03}, the discrete Poincar\'e inequalities for piecewise $H^1$ functions are proposed. However, the inverse, trace inequalities and the local quasi-uniformity for meshes under the shape-regular condition are used for the proof. Therefore, careful consideration of the results used in \cite{Bre03} may be necessary to remove the assumption that $\Omega$ is convex. 
	
\end{rem}

\begin{comment}
\section{Motivation}
We review classical dG methods to motivate our analysis. In this section only, we assume that a family of meshes $\{ \mathbb{T}_h\}$ has the (standard) shape-regular property. We define the calassical discontinuous $\mathbb{P}^1$ finite element space as
\begin{align*}
\displaystyle
V_{dc,h}^1 := \{ \varphi_h \in L^2(\Omega): \ v_h|_T \in \mathbb{P}^1(T) \ \forall T \in \mathbb{T}_h \}.
\end{align*}

The calassical dG methods are applied to many problems. We first focus on the symmetric interior penelty (SIP) method introduced by Arnord \cite{Arn82}, also see \cite{PieErn12,Riv08}. The SIP method is to find $u_h^{sip} \in V_{dc,h}^1$ such that

\end{comment}

\section{WOPSIP method for the Poisson equation} \label{sec=WOPSIP}
This section presents an analysis of the WOPSIP method for the Poisson equations on anisotropic meshes. In \cite{Ish24b}, we presented the error estimate in the energy norm $| \cdot |_{wop}$ in Section \ref{pena=norm} for the Stokes equation. We review the energy norm eror estimate and here is a new introduction to an error estimate in the $L^2$ norm.

\begin{comment}
We use the following inequality. Using the Cauchy--Schwarz and Jensen inequalities for $a_1,a_2,a_3 \in \mathbb{R}_+$ and $b_1 , b_2 \in \mathbb{R}_+ \cup \{ 0\}$,
\begin{align}
\displaystyle
&\sum_{i=1}^d a_1 (a_2+b_1)^{\frac{1}{2}}(a_3+b_2)^{\frac{1}{2}} \notag\\
&\leq c \left( \sum_{i=1}^d a_1^2 \right)^{\frac{1}{2}} \left\{ \left( \sum_{i=1}^d a_2^2 \right)^{\frac{1}{2}} + \left( \sum_{i=1}^d b_1^2 \right)^{\frac{1}{2}} \right\}^{\frac{1}{2}} \left\{ \left( \sum_{i=1}^d a_3^2 \right)^{\frac{1}{2}} + \left( \sum_{i=1}^d b_2^2 \right)^{\frac{1}{2}} \right\}^{\frac{1}{2}}. \label{newwop=1}
\end{align}

We introduce a Jensen-type inequality (see \cite[Exercise 12.1]{ErnGue21a}). Let $r,s$ be two nonnegative real numbers and $\{ x_i \}_{i \in I}$ be a finite sequence of nonnegative numbers. It then holds that
\begin{align}
\displaystyle
\begin{cases}
\left( \sum_{i \in I} x_i^s \right)^{\frac{1}{s}} \leq \left( \sum_{i \in I} x_i^r \right)^{\frac{1}{r}} \quad \text{if $r \leq s$},\\
\left( \sum_{i \in I} x_i^s \right)^{\frac{1}{s}} \leq \card(I)^{\frac{r-s}{rs}} \left( \sum_{i \in I} x_i^r \right)^{\frac{1}{r}} \quad \text{if $r \> s$}.
\end{cases} \label{jensen}
\end{align}
\end{comment}

\subsection{WOPSIP method}
We consider the WOPSIP method for the Poisson equation \eqref{poisson_eq} as follows. We aim to find $u_h^{wop} \in V_{dc,h}^{CR}$ such that
\begin{align}
\displaystyle
a_h^{wop}(u_h^{wop} , \varphi_h) = \ell_h(\varphi_h) \quad \forall \varphi_h \in V_{dc,h}^{CR}, \label{newwop=2}
\end{align}
where $\ell_h: V_{dc,h}^{CR} \to \mathbb{R}$ is defined in \eqref{lh=rem1}. Here, $a_h^{wop}: (V_{dc,h}^{CR} + H_0^1(\Omega) ) \times (V_{dc,h}^{CR} + H_0^1(\Omega) ) $ is defined as
\begin{align*}
\displaystyle
a_h^{wop}(v,w) &:= \int_{\Omega} \nabla_h v \cdot \nabla_h w dx + \sum_{F \in \mathcal{F}_h} \kappa_{F}  \int_F \Pi_F^{0} [\![ v]\!] \Pi_F^{0} [\![ w]\!] ds
\end{align*}
for all $(v,w) \in (V_{dc,h}^{CR} + H_0^1(\Omega) )  \times (V_{dc,h}^{CR} + H_0^1(\Omega) ) $. Recall that the parameter $\kappa_F$ is defined in \eqref{penalty0}. Using the H\"older's inequality, we obtain
\begin{align}
\displaystyle
|a_h^{wop}(v,w_h) | &\leq c |v|_{wop}  |w_h|_{wop} \quad \forall v \in V_{dc,h}^{CR} + H_0^1(\Omega) , \  \forall w_h \in V_{dc,h}^{CR}. \label{newwop=3}
\end{align}
As stated in Remark \ref{poin=rem1}, we impose that $\Omega$ is convex to obtain a stability estimate of the WOPSIP sheme. The stability estimate without convexity of the domain is still an open question.

\subsection{Energy norm error estimate}
The starting point for error analysis is the Second Strang Lemma, e.g. see \cite[Lemma 2.25]{ErnGue04}.

\begin{lem} \label{second=starang}
We assume that $\Omega$ is convex. Let $u \in H_0^1(\Omega)$ be the solution of \eqref{poisson_weak} and $u_h^{wop} \in V_{dc,h}^{CR}$ be the solution of \eqref{newwop=2}. It then holds that
\begin{align}
\displaystyle
|u - u_h^{wop}|_{wop}
&\leq \inf_{v_h \in V_{dc,h}^{CR}} |u - v_h|_{wop} + E_h(u), \label{newwop=4}
\end{align}
where
\begin{align}
\displaystyle
E_h(u) := \sup_{w_h \in V_{dc,h}^{CR}} \frac{|a_h^{wop}(u,w_h) - \ell_h(w_h) |}{|w_h|_{wop}}. \label{newwop=5}
\end{align}
\end{lem}

\begin{pf*}
Let $v_h \in V_{dc,h}^{CR}$. It holds that
\begin{align*}
\displaystyle
|v_h - u_h^{wop}|_{wop}^2
&= a_h^{wop}(v_h - u_h^{wop},v_h - u_h^{wop}) \\
&\hspace{-1cm} = a_h^{wop}(v_h - u,v_h - u_h^{wop}) + a_h^{wop}(u,v_h - u_h^{wop}) - \ell_h(v_h - u_h^{wop}) \\
&\hspace{-1cm} \leq c |u - v_h|_{wop}  |v_h - u_h^{wop}|_{wop} + | a_h^{wop}(u,v_h - u_h^{wop}) - \ell_h(v_h - u_h^{wop}) |,
\end{align*}
which leads to
\begin{align*}
\displaystyle
|v_h - u_h^{wop}|_{wop}
&\leq c |u - v_h|_{wop} + \frac{| a_h^{wop}(u,v_h - u_h^{wop}) - \ell_h(v_h - u_h^{wop}) |}{|v_h - u_h^{wop}|_{wop}} \\
&\leq  c |u - v_h|_{wop} + E_h(u).
\end{align*}
Then,
\begin{align*}
\displaystyle
|u - u_h^{wop}|_{wop}
&\leq |u - v_h|_{wop} + |v_h - u_h^{wop}|_{wop} \leq c  |u - v_h|_{wop}  + E_h(u).
\end{align*}
Hence, the target inequality \eqref{newwop=4} holds.
\qed
\end{pf*}

\begin{lem}[Best approximation] \label{best=approx}
We assume that $\Omega$ is convex. Let $u \in V_* := H_0^1(\Omega) \cap H^2(\Omega)$ be the solution of \eqref{poisson_weak}. Then,
\begin{align}
\displaystyle
\inf_{v_h \in V_{dc,h}^{CR}} |u - v_h|_{wop}
\leq c \left(  \sum_{i=1}^d  \sum_{T \in \mathbb{T}_h} h_i^2 \left \| \frac{\partial}{\partial r_i} \nabla u \right \|_{L^2(T)^d}^2 \right)^{\frac{1}{2}}. \label{newwop=6}
\end{align}
\end{lem}

\begin{pf*}
Let $v \in H_0^1(\Omega)$. From \cite[Theorem 7]{Ish24b},
\begin{align}
\displaystyle
\| \Pi_F^0 [\![  ({I}_{h}^{CR} v) - v] \!] \|^2_{L^2(F)} = 0. \label{newwop=7}
\end{align}
%Therefore, using the Jensen-type inequality \eqref{jensen}, \eqref{CR4} and \eqref{wop=6}, we obtain
Therefore, using \eqref{CR4} and \eqref{wop=6}, we obtain
\begin{align}
\displaystyle
 \inf_{v_h \in {V_{dc,h}^{CR}}}  |u-v_h|_{wop}
 &\leq  |u- {I}_{h}^{CR} u |_{wop}  = |u- {I}_{h}^{CR} u|_{H^1(\mathbb{T}_h)} \label{cr=jump}\\
 &\leq c \left(  \sum_{i=1}^d  \sum_{T \in \mathbb{T}_h} h_i^2 \left \| \frac{\partial}{\partial r_i} \nabla u \right \|_{L^2(T)^d}^2 \right)^{\frac{1}{2}}, \notag
% &\leq c \sum_{T \in \mathbb{T}_h} \sum_{i=1}^d {h}_i \left\| \frac{\partial }{\partial r_i} \nabla u \right \|_{L^2(T)^d},
\end{align}
which is the target inequality \eqref{newwop=6}.
\qed
\end{pf*}

The essential part for error estimates is the consistency error term \eqref{newwop=5}.

\begin{lem}[Asymptotic Consistency] \label{asy=con}
We assume that $\Omega$ is convex. Let $u \in V_* := H_0^1(\Omega) \cap H^2(\Omega)$ be the solution of \eqref{poisson_weak}. Let $\{ \mathbb{T}_h\}$ be a family of conformal meshes with the semi-regular property (Assumption \ref{neogeo=assume}). Let $T \in \mathbb{T}_h$ be the element with Conditions \ref{cond1} or \ref{cond2} and satisfy (Type \roman{sone}) in Section \ref{element=cond} when $d=3$. Then,
\begin{align}
\displaystyle
E_h(u)
&\leq c \left\{ \left( \sum_{i=1}^d \sum_{T \in \mathbb{T}_h} h_i^2 \left \| \frac{\partial}{\partial r_i} \nabla u \right \|_{L^2(T)^d}^2 \right)^{\frac{1}{2}} + h \| \varDelta u \|_{L^2(\Omega)} \right\} \notag \\
&\quad + c \left ( h |u|_{H^1(\Omega)} + h^{\frac{3}{2}} |u|_{H^1(\Omega)}^{\frac{1}{2}} \| \varDelta u\|_{L^2(\Omega)}^{\frac{1}{2}}  \right).   \label{newwop=8}
\end{align}
Furthermore, let $d=3$ and let {$T \in \mathbb{T}_h$} be the element with Condition \ref{cond2} and satisfy (Type \roman{stwo}) in Section \ref{element=cond}. It then holds that
\begin{align}
\displaystyle
E_h(u) 
&\leq c h \| \varDelta u \|_{L^2(\Omega)} + c \left ( h |u|_{H^1(\Omega)} + h^{\frac{3}{2}} |u|_{H^1(\Omega)}^{\frac{1}{2}} \| \varDelta u\|_{L^2(\Omega)}^{\frac{1}{2}}  \right). \label{newwop=9}
\end{align}
\end{lem}

\begin{pf*}
We first have
\begin{align*}
\displaystyle
 \div \mathcal{I}_{h}^{RT} \nabla u &= \Pi_h^{0} \div \nabla u =  \Pi_h^{0} \varDelta u.
\end{align*}
Let $w_h \in V_{dc,h}^{CR}$. Setting $w := \nabla u$ in \eqref{wop=3} yields
\begin{align*}
\displaystyle
&a_h^{wop}(u,w_h) - \ell_h(w_h) \\
&\quad =  \int_{\Omega} ( \nabla u- \mathcal{I}_{h}^{RT} \nabla u) \cdot \nabla_h w_{h} dx + \int_{\Omega} \left( \varDelta u -  \Pi_h^0 \varDelta u \right) w_{h}  dx  \\
&\quad \quad + \sum_{F \in \mathcal{F}_h^i} \int_{F}  \{ \! \{ \nabla u \} \!\}_{\omega,F} \cdot n_F \Pi_F^0 [\![ w_{h} ]\!]_F ds + \sum_{F \in \mathcal{F}_h^{\partial}} \int_{F} ( \nabla u \cdot n_F) \Pi_F^0 w_{h} ds  \\
&\quad =: I_1 + I_2 + I_3 + I_4. 
\end{align*}

Let $T \in \mathbb{T}_h$ be the element with Conditions \ref{cond1} or \ref{cond2} and satisfy (Type \roman{sone}) in Section \ref{element=cond} when $d=3$. Using the H\"older's inequality, the Cauchy--Schwarz inequality and the RT interpolation error \eqref{RT5}, the term $I_1$ is estimated as
\begin{align*}
\displaystyle
|I_1|
&\leq c  \sum_{T \in \mathbb{T}_h} \| \nabla u - \mathcal{I}_{h}^{RT} \nabla u \|_{L^2(T)} | w_{h} |_{H^1(T)} \\
&\leq c \sum_{T \in \mathbb{T}_h}\left(  \sum_{i=1}^d h_i \left \| \frac{\partial}{\partial r_i} \nabla u \right \|_{L^2(T)^d} + h_T \| \varDelta u \|_{L^2(T)} \right) | w_{h} |_{H^1(T)} \\
&\leq c \left\{ \left( \sum_{i=1}^d \sum_{T \in \mathbb{T}_h} h_i^2 \left \| \frac{\partial}{\partial r_i} \nabla u \right \|_{L^2(T)^d}^2 \right)^{\frac{1}{2}} + h \| \varDelta u \|_{L^2(\Omega)} \right\}  | w_{h} |_{wop}.
\end{align*}
Using the H\"older's inequality, the Cauchy--Schwarz inequality, the stability of $\Pi_h^0$ and the estimate \eqref{L2ortho}, the term $I_2$ is estimated as
\begin{align*}
\displaystyle
|I_2|
&= \left|  \int_{\Omega} \left( \varDelta u -  \Pi_h^0 \varDelta u \right) \left(  w_{h} - \Pi_h^0 w_{h}  \right) dx \right| \\
&\leq \sum_{T \in \mathbb{T}_h} \| \varDelta u -  \Pi_h^0 \varDelta u \|_{L^2(T)}  \| w_{h} - \Pi_h^0 w_{h} \|_{L^2(T)} \\
&\leq c h \| \varDelta u \|_{L^2(\Omega)}  | w_{h} |_{wop}.
\end{align*}
The triangle inequality, \eqref{wop=5} and \eqref{wop=6}, the terms $I_3$ and $I_4$ are estimated as
\begin{align*}
\displaystyle
|I_3| 
&\leq c |w_{h}|_{jwop} \Biggl(  h \|  \nabla u \|_{L^2(\Omega)^d}
+ h^{\frac{3}{2}} \| \nabla u \|_{L^2(\Omega)^d}^{\frac{1}{2}} |  \nabla u |_{H^1(\Omega)^d}^{\frac{1}{2}} \Biggr) \\
&\leq c |w_{h}|_{wop} \Biggl ( h |u|_{H^1(\Omega)} + h^{\frac{3}{2}} |u|_{H^1(\Omega)}^{\frac{1}{2}} \| \varDelta u\|_{L^2(\Omega)}^{\frac{1}{2}}  \Biggr), \\
|I_4|
&\leq c |w_{h}|_{wop} \Biggl ( h |u|_{H^1(\Omega)} + h^{\frac{3}{2}} |u|_{H^1(\Omega)}^{\frac{1}{2}} \| \varDelta u\|_{L^2(\Omega)}^{\frac{1}{2}}  \Biggr).
\end{align*}
Gathering the above inequalities yields the target inequality \eqref{newwop=8}.

Let $d=3$ and let {$T \in \mathbb{T}_h$} be the element with Condition \ref{cond2} and satisfy (Type \roman{stwo}) in Section \ref{element=cond}. Using the H\"older's inequality, the Cauchy--Schwarz inequality and the RT interpolation error \eqref{RT6}, the terms $I_1$ is estimated as
\begin{align*}
\displaystyle
|I_1|
&\leq c h  |  u |_{H^2(\Omega)}  | w_{h} |_{wop}, 
\end{align*}
which implies the target inequality \eqref{newwop=9} with \eqref{poisson=reg}.
\qed
\end{pf*}

From Lemmata \ref{second=starang}, \ref{best=approx} and \ref{asy=con}, we have the following energy norm error estimate.

\begin{thr} \label{energy=thr5}
We assume that $\Omega$ is convex. Let $u \in V_*$ be the solution of \eqref{poisson_weak}. Let $\{ \mathbb{T}_h\}$ be a family of conformal meshes with the semi-regular property (Assumption \ref{neogeo=assume}). Let $T \in \mathbb{T}_h$ be the element with Conditions \ref{cond1} or \ref{cond2} and satisfy (Type \roman{sone}) in Section \ref{element=cond} when $d=3$. Let $u_h^{wop} \in V_{dc,h}^{CR}$ be the solution of \eqref{newwop=2}. Then,
\begin{align}
\displaystyle
|u - u_h^{wop}|_{wop}
&\leq c \left\{ \left( \sum_{i=1}^d \sum_{T \in \mathbb{T}_h} h_i^2 \left \| \frac{\partial}{\partial r_i} \nabla u \right \|_{L^2(T)^d}^2 \right)^{\frac{1}{2}} + h \| \varDelta u \|_{L^2(\Omega)} \right\} \notag \\
&\quad + c \left ( h |u|_{H^1(\Omega)} + h^{\frac{3}{2}} |u|_{H^1(\Omega)}^{\frac{1}{2}} \| \varDelta u\|_{L^2(\Omega)}^{\frac{1}{2}}  \right).   \label{newwop=10}
\end{align}
Furthermore, let $d=3$ and let {$T \in \mathbb{T}_h$} be the element with Condition \ref{cond2} and satisfy (Type \roman{stwo}) in Section \ref{element=cond}. It then holds that
\begin{align}
\displaystyle
|u - u_h^{wop}|_{wop}
&\leq c h \| \varDelta u \|_{L^2(\Omega)} + c \left ( h |u|_{H^1(\Omega)} + h^{\frac{3}{2}} |u|_{H^1(\Omega)}^{\frac{1}{2}} \| \varDelta u\|_{L^2(\Omega)}^{\frac{1}{2}}  \right). \label{newwop=11}
\end{align}
\end{thr}

\subsection{$L^2$ norm error estimate}
This section presents the $L^2$ error estimate of the WOPSIP method.

\begin{thr}
We assume that $\Omega$ is convex. Let $u \in V_*$ be the solution of \eqref{poisson_weak}. Let $\{ \mathbb{T}_h\}$ be a family of conformal meshes with the semi-regular property (Assumption \ref{neogeo=assume}). Let $T \in \mathbb{T}_h$ be the element with Conditions \ref{cond1} or \ref{cond2} and satisfy (Type \roman{sone}) in Section \ref{element=cond} when $d=3$. Let $u_h^{wop} \in V_{dc,h}^{CR}$ be the solution of \eqref{newwop=2}. Then,
\begin{align}
\displaystyle
\| u - u_h^{wop} \|_{L^2(\Omega)}
&\leq c h \left\{ \left( \sum_{i=1}^d \sum_{T \in \mathbb{T}_h} h_i^2 \left \| \frac{\partial}{\partial r_i} \nabla u \right \|_{L^2(T)^d}^2 \right)^{\frac{1}{2}} + h \| \varDelta u \|_{L^2(\Omega)} \right\} \notag \\
&\quad + c h  \left ( h |u|_{H^1(\Omega)} + h^{\frac{3}{2}} |u|_{H^1(\Omega)}^{\frac{1}{2}} \| \varDelta u\|_{L^2(\Omega)}^{\frac{1}{2}}  \right). \label{newwop=est1}
\end{align}
Furthermore, let $d=3$ and let {$T \in \mathbb{T}_h$} be the element with Condition \ref{cond2} and satisfy (Type \roman{stwo}) in Section \ref{element=cond}. It then holds that
\begin{align}
\displaystyle
\| u - u_h^{wop} \|_{L^2(\Omega)}
&\leq c h^2  \| \varDelta u \|_{L^2(\Omega)} \notag \\
&\quad + c h  \left ( h |u|_{H^1(\Omega)} + h^{\frac{3}{2}} |u|_{H^1(\Omega)}^{\frac{1}{2}} \| \varDelta u\|_{L^2(\Omega)}^{\frac{1}{2}}  \right). \label{newwop=est2}
\end{align}
\end{thr}

\begin{pf*}
We set $e :=  u - u_h^{wop}$. Let $z \in V_*$ satisfy
\begin{align}
\displaystyle
a(\varphi , z) = \int_{\Omega} \varphi  e dx \quad \forall \varphi \in H_0^1(\Omega), \label{l2wop=13}
\end{align}
and $z_h^{wop} \in V_{dc,h}^{CR}$ satisfy 
\begin{align}
\displaystyle
a_h^{wop}(\varphi_h , z_h^{wop}) = \int_{\Omega} \varphi_h  e dx \quad \forall \varphi_h \in V_{dc,h}^{CR}. \label{l2wop=14}
% \int_{\Omega} \nabla_h \varphi_h \cdot \nabla_h z_h dx + \sum_{F \in \mathcal{F}_h} \kappa_{F}  \int_F \Pi_F^{0} [\![ \varphi_h]\!] \Pi_F^{0} [\![ z_h]\!] ds = \int_{\Omega} \varphi_h  e dx \quad \forall \varphi_h \in V_{dc,h}^{CR}. \label{l2wop=14}
\end{align}
Note that $|z|_{H^1(\Omega)} \leq c \| e \|_{L^2(\Omega)}$ and $|z|_{H^2(\Omega)} \leq \| e \|_{L^2(\Omega)}$. Then,
\begin{align*}
\displaystyle
\| e \|^2_{L^2(\Omega)}
&= \int_{\Omega} (u - u_h^{wop}) e dx = a(u,z) - a_h^{wop}(u_h^{wop},z_h^{wop}) \notag \\
&\hspace{-1cm} = a_{h}^{wop}(u -  u_h^{wop} , z - z_h^{wop} ) + a_{h}^{wop}( u -  u_h^{wop} ,  z_h^{wop}) + a_{h}^{wop}( u_h^{wop} ,  z - z_h^{wop} ) \notag \\
&\hspace{-1cm} = a_{h}^{wop}(u -  u_h^{wop} , z - z_h^{wop} ) \notag \\
&\hspace{-1cm} \quad + a_{h}^{wop}( u -  u_h^{wop} ,  z_h^{wop} - I_{h}^{CR} z) +   a_{h}^{wop}( u -  u_h^{wop} ,   I_{h}^{CR} z) \notag\\
&\hspace{-1cm} \quad + a_{h}^{wop}( u_h^{wop} - I_{h}^{CR} u ,  z - z_h^{wop} ) +a_{h}^{wop}( I_{h}^{CR} u ,  z - z_h^{wop} ) \\
&\hspace{-1cm} =: J_1 + J_2 + J_3 + J_4 + J_5.
\end{align*}
From Theorems \ref{DGCE=thr3} and \ref{energy=thr5}, we have
\begin{subequations} \label{l2wop=15}
\begin{align}
\displaystyle
 | z- z_h^{wop} |_{wop}
 &\leq c h \| e \|_{L^2(\Omega)}, \label{l2wop=15a} \\
 | z- I_{h}^{CR} z |_{H^1(\mathbb{T}_h)}
  &\leq c  h \| e \|_{L^2(\Omega)}, \label{l2wop=15b} \\
 \| z- I_{h}^{CR} z \|_{L^2(\Omega)}
  &\leq c  h^2 \| e \|_{L^2(\Omega)}. \label{l2wop=15c}
\end{align}
\end{subequations}

Let $T \in \mathbb{T}_h$ be the element with Conditions \ref{cond1} or \ref{cond2} and satisfy (Type \roman{sone}) in Section \ref{element=cond} when $d=3$.

Using \eqref{newwop=10} and \eqref{l2wop=15a}, $J_1$ can be estimated as
\begin{align*}
\displaystyle
|J_1| 
&\leq c  | u- u_h^{wop} |_{wop}  | z- z_h^{wop} |_{wop} \\
&\leq c h  \| e \|_{L^2(\Omega)} \left\{ \left( \sum_{i=1}^d \sum_{T \in \mathbb{T}_h} h_i^2 \left \| \frac{\partial}{\partial r_i} \nabla u \right \|_{L^2(T)^d}^2 \right)^{\frac{1}{2}} + h \| \varDelta u \|_{L^2(\Omega)} \right\} \notag \\
&\quad + c h  \| e \|_{L^2(\Omega)} \left ( h |u|_{H^1(\Omega)} + h^{\frac{3}{2}} |u|_{H^1(\Omega)}^{\frac{1}{2}} \| \varDelta u\|_{L^2(\Omega)}^{\frac{1}{2}}  \right).
\end{align*}
Using \eqref{cr=jump}, \eqref{newwop=10}, \eqref{l2wop=15a} and \eqref{l2wop=15b}, $J_2$ can be estimated as
\begin{align*}
\displaystyle
|J_2|
&\leq c \left \{ | u- u_h^{wop} |_{wop} \left(  | z_h^{wop} - z |_{wop} + | z - I_h^{CR} z|_{wop} \right) \right \} \\
&\leq c h  \| e \|_{L^2(\Omega)} \left\{ \left( \sum_{i=1}^d \sum_{T \in \mathbb{T}_h} h_i^2 \left \| \frac{\partial}{\partial r_i} \nabla u \right \|_{L^2(T)^d}^2 \right)^{\frac{1}{2}} + h \| \varDelta u \|_{L^2(\Omega)} \right\} \notag \\
&\quad + c h  \| e \|_{L^2(\Omega)} \left ( h |u|_{H^1(\Omega)} + h^{\frac{3}{2}} |u|_{H^1(\Omega)}^{\frac{1}{2}} \| \varDelta u\|_{L^2(\Omega)}^{\frac{1}{2}}  \right).
\end{align*}
By an analogous argument,
\begin{align*}
\displaystyle
|J_4|
%&\leq c \left \{ | z- z_h^{wop} |_{wop} \left(  | u_h^{wop} - u |_{wop} + | u - I_h^{CR} u|_{wop} \right) \right \} \\
&\leq c h  \| e \|_{L^2(\Omega)} \left\{ \left( \sum_{i=1}^d \sum_{T \in \mathbb{T}_h} h_i^2 \left \| \frac{\partial}{\partial r_i} \nabla u \right \|_{L^2(T)^d}^2 \right)^{\frac{1}{2}} + h \| \varDelta u \|_{L^2(\Omega)} \right\} \notag \\
&\quad + c h  \| e \|_{L^2(\Omega)} \left ( h |u|_{H^1(\Omega)} + h^{\frac{3}{2}} |u|_{H^1(\Omega)}^{\frac{1}{2}} \| \varDelta u\|_{L^2(\Omega)}^{\frac{1}{2}}  \right).
\end{align*}
Using the same argument with \eqref{wop=3} with $w := \nabla u$ and $\Pi_F^0 [\![  ({I}_{h}^{CR} z) - z] \!] = 0$,
\begin{align*}
\displaystyle
&\int_{\Omega} \left( \mathcal{I}_h^{RT} \nabla u \cdot \nabla_h (I_{h}^{CR} z - z ) + \Pi_h^0 \varDelta u  (I_{h}^{CR} z - z ) \right) dx = 0.
\end{align*}
Thus,
\begin{align*}
\displaystyle
J_3
&= a_{h}^{wop}( u  ,   I_{h}^{CR} z) - a_{h}^{wop}( u_h^{wop} ,   I_{h}^{CR} z) \\
&= \int_{\Omega} \nabla u \cdot \nabla_h ( I_{h}^{CR} z - z ) dx + \int_{\Omega} \varDelta u (I_{h}^{CR} z - z )dx \\
&=  \int_{\Omega} (\nabla u -  \mathcal{I}_h^{RT} \nabla u  ) \cdot \nabla_h ( I_{h}^{CR} z - z ) dx + \int_{\Omega} ( \varDelta u - \Pi_h^0 \varDelta u ) (I_{h}^{CR} z - z )dx.
\end{align*}
Using \eqref{RT5}, \eqref{l2wop=15b}, \eqref{l2wop=15c} and the stanility of the $L^2$-projection yields
\begin{align*}
\displaystyle
|J_3|
&\leq  c h  \| e \|_{L^2(\Omega)} \left\{ \left( \sum_{i=1}^d \sum_{T \in \mathbb{T}_h} h_i^2 \left \| \frac{\partial}{\partial r_i} \nabla u \right \|_{L^2(T)^d}^2 \right)^{\frac{1}{2}} + h \| \varDelta u \|_{L^2(\Omega)} \right\} \\
&\quad + c h^2  \| e \|_{L^2(\Omega)} \| \varDelta u \|_{L^2(\Omega)}.
\end{align*}
By an analogous argument,
\begin{align*}
\displaystyle
|J_5| 
&= \left| \int_{\Omega} (\nabla z -  \mathcal{I}_h^{RT} \nabla z  ) \cdot \nabla_h ( I_{h}^{CR} u - u ) dx + \int_{\Omega} ( \varDelta z - \Pi_h^0 \varDelta z ) (I_{h}^{CR} u - u )dx \right| \\
&\leq  c h  \| e \|_{L^2(\Omega)}  \left(  \sum_{i=1}^d  \sum_{T \in \mathbb{T}_h} h_i^2 \left \| \frac{\partial}{\partial r_i} \nabla u \right \|_{L^2(T)^d}^2 \right)^{\frac{1}{2}} + c h^2  \| e \|_{L^2(\Omega)} \| \varDelta u \|_{L^2(\Omega)}.
\end{align*}
Gathering the above inequalities yields the target inequality \eqref{newwop=est1}.

Let $d=3$ and let {$T \in \mathbb{T}_h$} be the element with Condition \ref{cond2} and satisfy (Type \roman{stwo}) in Section \ref{element=cond}. 

We use \eqref{newwop=11} instead of \eqref{newwop=10}, and \eqref{RT6} instead of \eqref{RT5} for estimates of $J_1 - J_5$:
\begin{align*}
\displaystyle
|J_1|
&\leq c h^2  \| e \|_{L^2(\Omega)} \| \varDelta u \|_{L^2(\Omega)} + c h  \| e \|_{L^2(\Omega)} \left ( h |u|_{H^1(\Omega)} + h^{\frac{3}{2}} |u|_{H^1(\Omega)}^{\frac{1}{2}} \| \varDelta u\|_{L^2(\Omega)}^{\frac{1}{2}}  \right), \\
|J_2|
&\leq c h^2  \| e \|_{L^2(\Omega)} \| \varDelta u \|_{L^2(\Omega)} + c h  \| e \|_{L^2(\Omega)} \left ( h |u|_{H^1(\Omega)} + h^{\frac{3}{2}} |u|_{H^1(\Omega)}^{\frac{1}{2}} \| \varDelta u\|_{L^2(\Omega)}^{\frac{1}{2}}  \right), \\
|J_3|
&\leq c h^2  \| e \|_{L^2(\Omega)} \| \varDelta u \|_{L^2(\Omega)},\\
|J_4|
&\leq c h^2  \| e \|_{L^2(\Omega)} \| \varDelta u \|_{L^2(\Omega)} + c h  \| e \|_{L^2(\Omega)} \left ( h |u|_{H^1(\Omega)} + h^{\frac{3}{2}} |u|_{H^1(\Omega)}^{\frac{1}{2}} \| \varDelta u\|_{L^2(\Omega)}^{\frac{1}{2}}  \right), \\
|J_5|
&\leq c h^2  \| e \|_{L^2(\Omega)} \| \varDelta u \|_{L^2(\Omega)}.
\end{align*}
Gathering the above inequalities yields the target inequality \eqref{newwop=est2}.
\qed
\end{pf*}

\section{Remarks on SIP methods} \label{rem=SIP}
%This section introduces remaks on SIP methods  for the Poisson equations on anisotropic meshes. To obtain a convergence analysis for the SIP method \eqref{rswip=1} without imposing the shape-regular condition for mesh partitions, as ususl, we estabilish consistency, boundedness for $a_h^{rsip}$, and discrete coercivity. However, we may not have an error estimate with an optimal order.

This section gives remaks on SIP and RSIP methods  for the Poisson equations without imposing the shape-regular condition for mesh partitions. We define the calassical discontinuous $\mathbb{P}^1$ finite element space as
\begin{align*}
\displaystyle
V_{dc,h}^1 &:= \{ \varphi_h \in L^2(\Omega): \ v_h|_T \in \mathbb{P}^1(T) \ \forall T \in \mathbb{T}_h \}, \\
V_{dc,h,*}^1 &:= V_* + V_{dc,h}^1,
\end{align*}
The SIP and RSIP methods are to find $u_h^{\natural} \in V_{dc,h}^{1}$ such that
\begin{align}
\displaystyle
a_h^{\natural}(u_h^{\natural} , \varphi_h) = \ell_h(\varphi_h) \quad \forall \varphi_h \in V_{dc,h}^{1}, \label{sip=1}
\end{align}
where $\ell_h: P_{dc,h}^{1} \to \mathbb{R}$ is defined in \eqref{lh=rem1} and $\natural \in \{sip,rsip \}$. Here, $a_h^{sip}: V_{dc,h,*}^{1} \times V_{dc,h}^{1}$ and $a_h^{rsip}: V_{dc,h,*}^{1} \times V_{dc,h}^{1}$ are defined as
\begin{align*}
\displaystyle
a_h^{sip}(v,w_h) &:= \int_{\Omega} \nabla_h v \cdot \nabla_h w_h dx + \sum_{F \in \mathcal{F}_h} \gamma^{sip} \kappa_{F*}  \int_F  [\![ v]\!]  [\![ w_{h}]\!] ds \\
&\quad - \sum_{F \in \mathcal{F}_h} \int_F \left( \{\!\{ \nabla_h v \}\!\}_{\omega} \cdot n_F [\![ w_h ]\!] + [\![ v ]\!] \{\!\{ \nabla_h w_h \}\!\}_{\omega} \cdot n_F \right ) ds, \\
a_h^{rsip}(v,w_h) &:= \int_{\Omega} \nabla_h v \cdot \nabla_h w_h dx + \sum_{F \in \mathcal{F}_h} \gamma^{rsip} \kappa_{F*}  \int_F \Pi_F^{0} [\![ v]\!] \Pi_F^{0} [\![ w_{h}]\!] ds \\
&\quad - \sum_{F \in \mathcal{F}_h} \int_F \left( \{\!\{ \nabla_h v \}\!\}_{\omega} \cdot n_F [\![ w_h ]\!] + [\![ v ]\!] \{\!\{ \nabla_h w_h \}\!\}_{\omega} \cdot n_F \right ) ds
\end{align*}
for all $(v,w_h) \in V_{dc,h,*}^{1} \times V_{dc,h}^{1}$. The penalty parameters $\gamma^{sip}$ and $\gamma^{rsip}$ are large enough. One generally estabilish consistency, boundedness for bilinear forms, and discrete coercivity to obtain a convergence analysis for the SIP and RSIP methods, e.g. see \cite[Chapter 4]{PieErn12}. However, we may not get an optimal order of error estimates for the consistenct term,
\begin{align*}
\displaystyle
 \sum_{F \in \mathcal{F}_h} \int_F  \{\!\{ \nabla_h v \}\!\}_{\omega} \cdot n_F [\![ w_h ]\!] ds \quad \forall (v,w_h) \in V_{dc,h,*}^{1} \times V_{dc,h}^{1}.
\end{align*}
Suppose that $F \in \mathcal{F}_h^i$ with $F = T_1 \cap T_2$, $T_1,T_2 \in \mathbb{T}_h$. Using the trace (Lemma \ref{lem=trace}) and the H\"older inequalities, the weighted average gives the following estimate for the term $\int_F \{\! \{ \nabla_h v \}\! \}_{\omega} \cdot n_F [\![ w_h ]\!] ds$. 
\begin{align}
\displaystyle
&\left|  \int_F  \{\!\{ \nabla_h v \}\!\}_{\omega} \cdot n_F [\![ w_h ]\!] ds  \right| \notag \\
&\quad \leq c \Bigg \{ \left(  \| \nabla v|_{T_1} \|_{L^2(T_1)^d} + h_{T_1}^{\frac{1}{2}}  \| \nabla v|_{T_1} \|_{L^2(T_1)^d}^{\frac{1}{2}} | \nabla v|_{T_1} |_{H^1(T_1)^d}^{\frac{1}{2}} \right)^2 \notag \\
&\quad \quad  + \left(  \| \nabla v|_{T_2} \|_{L^2(T_2)^d} + h_{T_2}^{\frac{1}{2}}  \| \nabla v|_{T_2} \|_{L^2(T_2)^d}^{\frac{1}{2}} | \nabla v|_{T_2} |_{H^1(T_2)^d}^{\frac{1}{2}}  \right)^{2} \Biggr\}^{\frac{1}{2}} \notag \\
&\quad \quad \times \left( \omega_{T_1,F}^2 { \ell_{T_1,F}}^{-1} + \omega_{T_2,F}^2 {\ell_{T_2,F}}^{-1} \right)^{\frac{1}{2}} \|  [\![ w_{h} ]\!] \|_{L^2(F)}, \label{sip=2}
\end{align}
where $\omega_{T_1,F}$ and $\omega_{T_2,F}$ are defined in \eqref{weight}, and $\ell_{T_1,F}$ and $\ell_{T_2,F}$ are defined in the inequality \eqref{trace=vec}. Let $I_T^L: \mathcal{C}(T) \to \mathbb{P}^1(T)$ for any $T \in \mathbb{T}_h$ be the usual Lagrange interpolation opertor,  and let $u \in V_*$ the exact solution of \eqref{poisson_eq}. We set $v := u - I_T^L u$. When $d=3$, due to \cite[Corolllary 1]{IshKobTsu21c}, ubder the semi-regular condition \eqref{NewGeo}, 
\begin{align*}
\displaystyle
|u - I_T^L u|_{H^1(T)}
&\leq c \sum_{i=1}^3 h_i \left | \frac{\partial u}{\partial r_i} \right |_{H^{1}(T)},
\end{align*}
which leads to
\begin{align*}
\displaystyle
 &\| \nabla (u - I_T^L u) \|_{L^2(T)^3} + h_{T}^{\frac{1}{2}}  \| \nabla (u - I_T^L u) \|_{L^2(T)^3}^{\frac{1}{2}} | \nabla (u - I_T^L u)|_{H^1(T)^3}^{\frac{1}{2}}  \\
 &\quad \leq c \sum_{i=1}^3 h_i \left | \frac{\partial u}{\partial r_i} \right |_{H^{1}(T)} + c h_T^{\frac{1}{2}} \left( \sum_{i=1}^3 h_i \left | \frac{\partial u}{\partial r_i} \right |_{H^{1}(T)} \right)^{\frac{1}{2}} |u|_{H^2(T)}^{\frac{1}{2}} \\
 &\quad \leq c h_T^{\frac{1}{2}} \left( \sum_{i=1}^3 h_i \left | \frac{\partial u}{\partial r_i} \right |_{H^{1}(T)} \right)^{\frac{1}{2}} |u|_{H^2(T)}^{\frac{1}{2}} \leq c h_T  |u|_{H^2(T)}.
\end{align*}
Therefore, even if anisotropic meshes are used, the computational efficiency may not increase much for the SIP method.

In case the RSIP method, the estimate on consistency term is more complicated. For all $(v,w_h) \in V_{dc,h,*}^{1} \times V_{dc,h}^{1}$,
\begin{align}
\displaystyle
&\int_F  \{\!\{ \nabla_h v \}\!\}_{\omega} \cdot n_F [\![ w_h ]\!] ds \notag\\
&= \int_F  \{\!\{ \nabla_h v \}\!\}_{\omega} \cdot n_F \Pi_F^0 [\![ w_h ]\!] ds
+ \int_F  \{\!\{ \nabla_h v \}\!\}_{\omega} \cdot n_F \left( [\![ w_h ]\!] - \Pi_F^0 [\![ w_h ]\!] \right) ds. \label{sip=3}
\end{align}
The first term of the R.H.S in \eqref{sip=3} can be estimated as \eqref{sip=2}, where the CR interpolation operator $I_T^{CR}$ is used instead of $I_T^L$. Therefore, even if anisotropic meshes are used, the computational efficiency for the RSIP method also may not increase much. In an estimate of the secoond term of the R.H.S in \eqref{sip=3}, we use scaling argument on anisotropic meshes, c.f. \cite[Lemma 4]{IshKobTsu21b}.

\section{Numerical experiments} \label{numerical=sec}
Under construction.

\begin{comment}

\begin{align*}
\displaystyle
\end{align*}

\begin{align}
\displaystyle
\end{align}

\begin{subequations}
\begin{align}
\displaystyle
\end{align}
\end{subequations}

\begin{flushright} 
$\square$
\end{flushright}

\begin{pf*}
\begin{flushright} 
$\square$
\end{flushright}
\end{pf*}

\begin{nnote}
\begin{flushright} 
$\square$
\end{flushright}
\end{nnote}

\begin{description}
  \item[(\Roman{lone})] 
  \item[(\Roman{ltwo})] 
  \item[(\Roman{lthree})] 
\end{description}

\begin{itemize}
 \item 
 \item 
 \item
\end{itemize}

\begin{enumerate}
 \item 
 \item 
 \item 
\end{enumerate}

\end{comment}

%
%**** Reference ****
%

\end{document}